# UNIQUENESS OF SOLUTIONS OF THE STOCHASTIC NAVIER–STOKES EQUATION WITH INVARIANT MEASURE GIVEN BY THE ENSTROPHY

By S. Albeverio[1] and B. Ferrario[2]

*Universität Bonn and Università di Pavia*

A stochastic Navier–Stokes equation with space-time Gaussian white noise is considered, having as infinitesimal invariant measure a Gaussian measure $\mu_\nu$ whose covariance is given in terms of the enstrophy. Pathwise uniqueness for $\mu_\nu$-a.e. initial velocity is proven for solutions having $\mu_\nu$ as invariant measure.

**1. Introduction.** We are interested in the stochastic Navier–Stokes equation with a space-time white noise. We consider the spatial domain to be the torus $\mathbb{T}^2 = [0, 2\pi]^2$ (hence periodic boundary conditions are assumed). In [1] it has been shown that there exists an infinitesimal invariant measure associated to this stochastic equation; this is a Gaussian measure $\mu_\nu$, with covariance given in terms of the enstrophy (and of the viscosity parameter $\nu$). Existence of a solution has been proven in two different ways: [1] considers a weak solution and [9] a strong solution (weak and strong are to be understood in the probabilistic sense). The common point of these papers is that the solution is obtained as the limit of Galerkin approximations. No result of uniqueness has been given in [1], whereas [9] shows existence and uniqueness in a smaller class than the natural one to consider for this problem. Indeed, the statement of Theorem 5.1 in [9] involves an auxiliary process (denoted by $z$ in Section 4), not appearing in the given stochastic Navier–Stokes equation, and for this reason the definition of uniqueness given in

Received November 2002; revised May 2003.

[1]Also affiliated with SFB 611, Bonn; BiBoS, Bielefeld; Dipartimento di Matematica, Università di Trento, I-38050 Povo; CERFIM, Locarno; Accademia di Architettura, USI, CH-6850 Mendrisio.

[2]Supported by the Alexander von Shumboldt Stiftung.

[9] is not the natural definition to consider and does not coincide with the pathwise uniqueness we prove in the present paper, as we will explain in Section 4.

The aim of this paper is to prove uniqueness of the solutions of this stochastic Navier–Stokes equation with a space-time white noise, in the same class where existence holds. Precisely, we will deal with processes with $\mathbb{P}$-a.e. path $u \in C([0,\infty); \mathcal{B}_{pq}^{-s})$ (with $\mathcal{B}_{pq}^{-s}$ being a certain Besov space specified below) and having $\mu_\nu$ as invariant measure.

Finally, we want to remark that the stochastic Navier–Stokes equation in a two-dimensional domain and with space-time Gaussian white noise has been discussed in some papers in the last years. Anyway, the only expression known for an invariant measure is that of the centered Gaussian measure $\mu_\nu$ considered in this paper too. No other invariant measures are known with this space-time Gaussian white noise. However, the (deterministic) 2D-Euler equation has many invariant measures, including all the measures $\mu_\nu$ (for $\nu > 0$) (see, e.g., [3] for a review on these invariant measures).

As to the structure of this paper, in Section 2 we shall introduce the two-dimensional Navier–Stokes equation and define the mathematical setting. In Section 3 the Gaussian measure $\mu_\nu$ of the enstrophy will be defined and the main properties of the nonlinear operator $B$ with respect to $\mu_\nu$ will be presented. The uniqueness result will be proven in Section 4. Two results used in the proofs will be given in the Appendix.

## 2. The Navier–Stokes equation.

We consider the equations governing the motion of a homogeneous incompressible viscous fluid in the two-dimensional torus

$$\frac{\partial}{\partial t}u(t,\xi) - \nu\Delta u(t,\xi) + [u(t,\xi) \cdot \nabla]u(t,\xi) - \nabla p(t,\xi) = f(t,\xi),$$

(2.1)  $\nabla \cdot u(t,\xi) = 0,$

$$u(0,\xi) = x(\xi),$$

with periodic boundary condition. The definition domains of the variables are $t \geq 0, \xi \in \mathbb{T}^2$. The unknowns are the velocity vector field $u = u(t,\xi)$ and the scalar pressure field $p = p(t,\xi)$. Here $\Delta = \frac{\partial^2}{\partial \xi_1^2} + \frac{\partial^2}{\partial \xi_2^2}$, $\nabla = (\frac{\partial}{\partial \xi_1}, \frac{\partial}{\partial \xi_2})$, $\xi = (\xi_1, \xi_2)$ and "$\cdot$" is the scalar product in $\mathbb{R}^2$. The viscosity $\nu$ is a strictly positive constant; $x$ and $f$ are the data.

We define the mathematical setting as follows. Consider any periodic divergence-free vector distribution $u$. Since $\nabla \cdot u = 0$, there exists a periodic scalar distribution $\psi$, called the stream function, such that

$$u = \nabla^\perp \psi \equiv \left(-\frac{\partial \psi}{\partial \xi_2}, \frac{\partial \psi}{\partial \xi_1}\right).$$

(2.2)

Decomposing $\psi$ in Fourier series with respect to the complete orthonormal system in $L_2(\mathbb{T}^2)$ given by $\{\frac{1}{2\pi}e^{ik\cdot\xi}\}_{k\in\mathbb{Z}^2}$

$$\psi(\xi) = \sum_{k\in\mathbb{Z}^2} \psi_k \frac{e^{ik\cdot\xi}}{2\pi}, \qquad \psi_k \in \mathbb{C}, \ \overline{\psi}_k = \psi_{-k},$$

by (2.2) we get that $u$ has the following Fourier series representation:

$$(2.3) \qquad u(\xi) = \sum_{k\in\mathbb{Z}_0^2} u_k e_k(\xi), \qquad u_k \in \mathbb{C}, \ \overline{u}_k = u_{-k},$$

where $e_k(\xi) = \frac{k^{\perp}}{2\pi|k|}e^{ik\cdot\xi}$. Here $k^{\perp} = (-k_2, k_1)$, $|k| = \sqrt{k_1^2 + k_2^2}$ and $\mathbb{Z}_0^2 = \{k \in \mathbb{Z}^2 : |k| \neq 0\}$. We define also $\mathbb{Z}_+^2 = \{k \in \mathbb{Z}_0^2 : k_1 > 0 \text{ or } \{k_1 = 0, k_2 > 0\}\}$.

Note that $\{e_k\}_{k\in\mathbb{Z}_0^2}$ is a complete orthonormal system of the eigenfunctions (with corresponding eigenvalues $|k|^2$) of the operator $-\Delta$ in $[L_2^{\mathrm{div}}(\mathbb{T}^2)]^2 = \{u \in [L_2(\mathbb{T}^2)]^2 : \nabla \cdot u = 0$, with the normal component of $u$ being periodic on $\partial\mathbb{T}^2\}$.

Each $e_k$ is a periodic divergence-free $C^{\infty}$-vector function. The convergence of the series (2.3) depends on the regularity of the vector function $u$, and can be used to define Sobolev spaces as in the following definition.

Let $\mathcal{U}'$ be the space of zero mean value periodic divergence-free vector distributions. Any element $u \in \mathcal{U}'$ is uniquely defined by the sequence of the coefficients $\{u_k\}_{k\in\mathbb{Z}_+^2}$; indeed, by duality, $u_k = \langle u, e_{-k}\rangle$, since each $e_k$ is a periodic divergence-free and infinitely differentiable function. Following [5], we define the periodic divergence-free vector Sobolev spaces, $s \in \mathbb{R}, 1 \leq p \leq \infty$,

$$\mathcal{H}_p^s = \left\{ u = \sum_{k\in\mathbb{Z}_0^2} u_k e_k \in \mathcal{U}' : \sum_k u_k |k|^s e_k(\cdot) \in L_p(\mathbb{T}^2) \right\}$$

and the periodic divergence-free Besov spaces as real interpolation spaces

$$\mathcal{B}_{pq}^s = (\mathcal{H}_p^{s_0}, \mathcal{H}_p^{s_1})_{\theta,q}, \qquad s \in \mathbb{R}, \ 1 \leq p, q \leq \infty,$$
$$s = (1-\theta)s_0 + \theta s_1, \qquad 0 < \theta < 1.$$

In particular, $\mathcal{B}_{22}^s = \mathcal{H}_2^s$. (For the theory of interpolation spaces see, e.g., [5].) Moreover, $\mathcal{U}' = \bigcup_{s\in\mathbb{R}, 1\leq p\leq\infty} \mathcal{H}_p^s$ with the inductive topology.

$\{e_k\}_{k\in\mathbb{Z}_0^2}$ is a complete orthonormal system in the space $\mathcal{H}_2^0$. It follows that the Hilbert space $\mathcal{H}_2^s$ is isomorphic to the space of complex valued sequences $\{u_k\}_{k\in\mathbb{Z}_0^2}$ such that $\sum_k |u_k|^2 |k|^{2s} < \infty$.

We define the Stokes operator as

$$A = -\Delta,$$

which is a linear operator in $\mathcal{H}_p^s$ with domain $\mathcal{H}_p^{s+2}$. It is an isomorphism from $\mathcal{H}_p^{s+2}$ to $\mathcal{H}_p^s$, $s \in \mathbb{R}$, $1 \le p < \infty$. For $u = \sum_k u_k e_k$, we have $Au = \sum_k u_k |k|^2 e_k$. Let $\Pi$ be the projector operator from the space of periodic vectors onto the space of periodic divergence-free vectors. Applying $\Pi$ to both sides of the first equation in the Navier–Stokes system, we get rid of the pressure term. The bilinear operator $B$ is defined as

$$
\begin{aligned}
B(u, v) &= \Pi\left[(u \cdot \nabla)v\right] \\
&= \Pi[\nabla \cdot (u \otimes v)] \qquad \text{(by the divergence-free condition)} \\
&= \Pi\left[\begin{pmatrix} \partial_1 \\ \partial_2 \end{pmatrix} \cdot \begin{pmatrix} u_1 v_1 & u_1 v_2 \\ u_2 v_1 & u_2 v_2 \end{pmatrix}\right]
\end{aligned}
$$

whenever it makes sense. For instance, a classical result is that $B : \mathcal{H}_2^1 \times \mathcal{H}_2^1 \to \mathcal{H}_2^{-1}$ (see, e.g., [13]). The (optimal) regularity of $B$ is the key point to solve the Navier–Stokes equation, both in the deterministic and in the stochastic case.

For less regular vectors $u$ and $v$, estimates on $B$ are given in Besov spaces (see, e.g., [6, 7]). This is useful in solving the stochastic Navier–Stokes equation with space-time white noise, as shown in [9].

We shall very often write $B(u)$ for the quadratic term $B(u, u)$.

The stochastic Navier–Stokes equation in which we are interested has the following abstract Itô form:

$$
\begin{aligned}
(2.4) \qquad & du(t) + [\nu Au(t) + B(u(t))]\, dt = dw(t), \qquad t > 0, \\
& u(0) = x.
\end{aligned}
$$

$\{w(t)\}_{t \ge 0}$ is a Wiener process, defined on a complete probability space $(\Omega, \mathcal{F}, \mathbb{P})$ with filtration $\{\mathcal{F}_t\}_{t \ge 0}$, which is cylindric in the space of finite energy $\mathcal{H}_2^0$; that is,

$$
w(t) = \sum_{k \in \mathbb{Z}_0^2} \beta_k(t) e_k,
$$

where $\{\beta_k\}_{k \in \mathbb{Z}_0^2}$ is a sequence of standard independent complex valued Wiener processes with $\beta_{-k} = \overline{\beta}_k$. This is a process with continuous paths taking values in $\mathcal{H}_2^\sigma$ for any $\sigma < -1$ (see, e.g., [10]). In other terms, $dw(t)$ is a Gaussian space-time white noise. We shall denote by $\mathbb{E}$ the expectation with respect to the measure $\mathbb{P}$.

The equation for the Fourier components is obtained by multiplying the first equation (2.4) by $e_{-k}(\xi)$ and integrating over the torus $\mathbb{T}^2$. We obtain, for any $k \in \mathbb{Z}_0^2$,

$$
\begin{aligned}
& du_k(t) + [\nu |k|^2 u_k(t) + B_k(u(t))]\, dt = d\beta_k(t), \qquad t > 0, \\
& u_k(0) = x_k,
\end{aligned}
$$

where

$$B_k(u) = \sum_{h \in \mathbb{Z}_0^2, h \neq k} c_{h,k} u_h u_{k-h},$$

$$c_{h,k} = -\frac{1}{4\pi} \frac{(h^\perp \cdot k)}{|h||k-h|} \left[ |k| - 2\frac{(h \cdot k)}{|k|} \right].$$

**3. The Gaussian invariant measure given by the enstrophy (and viscosity parameter).** We shall consider a certain centered Gaussian measure $\mu_\nu$ on the space $\mathcal{U}'$ of complex valued sequences $\{u_k\}_{k \in \mathbb{Z}_+^2}$. $\mu_\nu$ is heuristically defined as the infinite product of (complex valued) centered Gaussian measures

$$(3.1) \qquad d\mu_\nu(u) = \frac{1}{Z} \underset{k \in \mathbb{Z}_+^2}{\times} e^{-2\nu|k|^2|u_k|^2} \, du_k$$

($|u_k|^2 = x_k^2 + y_k^2, du_k = dx_k \, dy_k$ for $u_k = x_k + iy_k$, $x_k, y_k \in \mathbb{R}$; $Z$ is a normalization factor).

Rigorously, $\mu_\nu$ is the mean zero Gaussian measure having as covariance the scalar product $(u, v)_\nu = \frac{1}{2\nu} \sum_{k \in \mathbb{Z}_+^2} |k|^{-2} u_k \overline{v_k}$. In particular,

$$\mathbb{E}_{\mu_\nu}[u_k \overline{u_j}] = \begin{cases} \dfrac{1}{2\nu|k|^2}, & \text{if } k = j, \\ 0, & \text{if } k \neq j. \end{cases}$$

The quantity in the exponent of the heuristic Gaussian density in (3.1) is the enstrophy $\mathcal{S}$ associated to the velocity field $u$: $\mathcal{S}(u) = \int_{\mathbb{T}^2} |\nabla^\perp \cdot u(\xi)|^2 \, d\xi \equiv 2 \sum_{k \in \mathbb{Z}_+^2} |k|^2 |u_k|^2$. In this sense, $\mu_\nu$ is the Gaussian measure given in terms of the enstrophy (and of the viscosity parameter $\nu$).

Let us characterize the support of the measure $\mu_\nu$. We have, for any integer $n$,

$$\mathbb{E}_{\mu_\nu}\left( \|u\|^{2n}_{\mathcal{H}_{2n}^{-s}} \right) = \int_{\mathcal{U}'} \|u\|^{2n}_{\mathcal{H}_{2n}^{-s}} \, d\mu_\nu(u)$$

$$= \int_{\mathcal{U}'} \left\| \sum_k u_k e_k \right\|^{2n}_{\mathcal{H}_{2n}^{-s}} \, d\mu_\nu(u)$$

$$= \int_{\mathcal{U}'} \left( \int_{\mathbb{T}^2} \left| \sum_k u_k |k|^{-s} e_k(\xi) \right|^{2n} d\xi \right) d\mu_\nu(u)$$

$$(3.2) \qquad = \int_{\mathbb{T}^2} \left( \int_{\mathcal{U}'} \left| \sum_k u_k |k|^{-s} e_k(\xi) \right|^{2n} d\mu_\nu(u) \right) d\xi$$

$$= c_n \int_{\mathbb{T}^2} \left[ \sum_k |k|^{-2s} |e_k(\xi)|^2 (\mathbb{E}_{\mu_\nu} |u_k|^2) \right]^n d\xi$$

$$= c_n' \left[ \sum_k |k|^{-2s} \mathbb{E}_{\mu_\nu} |u_k|^2 \right]^n$$

for some constants $c_n, c_n' > 0$. In these calculations we have used that, for any $\gamma_k \in \mathbb{C}$,

$$(3.3) \qquad \mathbb{E}_{\mu_\nu} \left| \sum_k u_k \gamma_k \right|^{2n} = \frac{(2n)!}{2^n n!} \left[ \sum_k |\gamma_k|^2 \mathbb{E}_{\mu_\nu}(|u_k|^2) \right]^n$$

and the fact that $|e_k(\xi)| = \frac{1}{2\pi}$ for any $\xi \in \mathbb{T}^2$.

Since $\mathbb{E}_{\mu_\nu}(|u_k|^2) = \frac{1}{2\nu|k|^2}$, the above calculation implies that there exists a positive constant $c_n''$ such that

$$\mathbb{E}_{\mu_\nu} \left( \|u\|_{\mathcal{H}_{2n}^{-s}}^2 \right) \leq \left( \mathbb{E}_{\mu_\nu} \|u\|_{\mathcal{H}_{2n}^{-s}}^{2n} \right)^{1/n} \leq \frac{1}{\nu} c_n'' \sum_{k \in \mathbb{Z}_0^2} \frac{1}{|k|^{2+2s}}.$$

The latter series converges as soon as $s > 0$. Hence $\mu_\nu(\mathcal{H}_{2n}^{-s}) = 1$ for any $s > 0$ and integer $n$. Since we are in a bounded spatial domain, we have the embedding $\mathcal{H}_{2(n+1)}^{-s} \subset \mathcal{H}_q^{-s} \subset \mathcal{H}_{2n}^{-s}$ for $2n < q < 2(n+1)$. Therefore,

$$\mu_\nu(\mathcal{H}_q^{-s}) = 1 \qquad \forall s > 0, \ 1 \leq q < \infty.$$

We remark that it was already known that the space $\mathcal{H}_2^0$ of finite-energy velocity vectors does not have full measure with respect to $\mu_\nu$; in fact, one has even $\mu_\nu(\mathcal{H}_2^0) = 0$ (see [4]).

We want to get Besov spaces of full measure $\mu_\nu$. First, we have the embedding

$$\mathcal{H}_q^{-s} \subseteq \mathcal{B}_{qq}^{-s}, \qquad 2 \leq q < \infty$$

(see [5], Theorem 6.4.4). Hence $\mu_\nu(\mathcal{B}_{qq}^{-s}) = 1$ for any $s > 0$ and $2 \leq q < \infty$. Moreover,

$$\mathcal{B}_{22}^{-s} \subset \mathcal{B}_{2q}^{-s}, \qquad 2 < q \leq \infty$$

(see [5], Theorem 6.2.4). Hence $\mu_\nu(\mathcal{B}_{2q}^{-s}) = 1$ for any $s > 0$ and $2 \leq q < \infty$. By interpolation, for $0 < \theta < 1$,

$$(\mathcal{B}_{qq}^{-s_0}, \mathcal{B}_{2q}^{-s_1})_{[\theta]} = \mathcal{B}_{pq}^{-s}$$

with $-s = (1-\theta)(-s_0) + \theta(-s_1)$ and $\frac{1}{p} = \frac{1-\theta}{q} + \frac{\theta}{2}$ (see [5], Theorem 6.4.5). This implies that, given $q$, for any $s > 0$ there exist $\theta \in (0,1)$ and $s_0, s_1 > 0$ such that the above interpolation holds. Necessarily we have $2 < p < q$. Hence $\mathcal{B}_{qq}^{-s_0} \cap \mathcal{B}_{2q}^{-s_1} \subset \mathcal{B}_{pq}^{-s}$, giving

$$(3.4) \qquad \mu_\nu(\mathcal{B}_{pq}^{-s}) = 1 \qquad \forall\, s > 0,\ 2 < p < q < \infty.$$

Summing up, we have proven the following result.

PROPOSITION 3.1. *For any viscosity $\nu > 0$,*

$$\mu_\nu(\mathcal{B}_{pq}^{-s}) = 1 \qquad \forall\, s > 0,\ 2 \le p \le q < \infty.$$

REMARK 3.1. With calculation similar to (3.2), we can obtain that, $\mathbb{P}$-a.s., the paths of the Wiener process $w(t) \in \mathcal{H}_p^{-1-s}$ for $s > 0$ and $1 \le p < \infty$.

We present now an estimate of the quadratic term $B$, useful in the following.

PROPOSITION 3.2. *For any viscosity $\nu > 0$, we have*

$$(3.5) \qquad \int \|B(u)\|_{\mathcal{H}_2^{-r-1}}^\rho \, d\mu_\nu(u) < \infty \qquad \forall\, r > 0,\ 1 \le \rho < \infty.$$

PROOF. Let us start by considering the case $\rho = 2$. We have that $B(u) = \sum_k B_k(u) e_k$ is defined as the limit in $\mathcal{H}_2^{-r-1}$ of $B^N(u) := \sum_{|k| \le N} B_k^N(u) e_k$, with $B_k^N(u) = \sum_{|h|, |k-h|, |k| \le N} c_{h,k} u_h u_{k-h}$. It will be shown that this limit exists in $\mathcal{L}^2(\mu_\nu)$ and that $\|B(u)\|_{\mathcal{H}_2^{-r-1}}^2 = \sum_k |B_k(u)|^2 |k|^{2(-r-1)}$. Let us compute the following integral with respect to the measure $\mu_\nu$:

$$\int \|B(u)\|_{\mathcal{H}_2^{-r-1}}^2 \, d\mu_\nu(u)$$

$$= \int \sum_{k \in \mathbb{Z}_0^2} |k|^{2(-r-1)} \left| \sum_{h \in \mathbb{Z}_0^2, h \ne k} c_{h,k} u_h u_{k-h} \right|^2 d\mu_\nu(u)$$

$$= \sum_{k \in \mathbb{Z}_0^2} |k|^{2(-r-1)} \int \sum_{h,h'} c_{h,k} c_{h',k} u_h u_{k-h} \overline{u_{h'} u_{k-h'}} \, d\mu_\nu(u)$$

$$= \sum_{k \in \mathbb{Z}_0^2} |k|^{2(-r-1)} \sum_{h \in \mathbb{Z}_0^2, h \ne k} (c_{h,k}^2 + c_{h,k} c_{k-h,k}) \frac{1}{2\nu |h|^2} \frac{1}{2\nu |k-h|^2},$$

where we have used the Fubini–Tonelli theorem to interchange the summations over $k$ and over $h$ with the integral.

Let us notice that the coefficients $c_{h,k}$ are such that $c_{h,k} = c_{k-h,k}$ and

$$
\begin{aligned}
c_{h,k}^2 &= \frac{1}{(4\pi)^2} \frac{|h^\perp \cdot k|^2}{|h|^2|k-h|^2|k|^2}[(k-h)\cdot k - h\cdot k]^2 \\
&\leq \frac{2}{(4\pi)^2} \frac{|h^\perp \cdot k|^2}{|h|^2|k-h|^2|k|^2}|(k-h)\cdot k|^2 + \frac{2}{(4\pi)^2} \frac{|(k-h)^\perp \cdot k|^2}{|h|^2|k-h|^2|k|^2}|h\cdot k|^2 \\
&\leq \frac{2}{(4\pi)^2} 2|k|^2.
\end{aligned}
$$

Then, continuing the estimates on the quadratic term, we get

$$
\begin{aligned}
\int \|B(u)\|_{\mathcal{H}_2^{-r-1}}^2 \, d\mu_\nu(u) &\leq \frac{1}{8\pi^2\nu^2} \sum_{k,h\in\mathbb{Z}_0^2, k\neq h} \frac{1}{|k|^{2r}|h|^2|k-h|^2} \\
&\leq \frac{1}{8\pi^2\nu^2} \sum_{k\in\mathbb{Z}_0^2} \frac{1}{|k|^{2r}} \sum_{h\in\mathbb{Z}_0^2, h\neq k} \frac{1}{|h|^2|k-h|^2} \\
&\leq \frac{c}{8\pi^2\nu^2} \sum_{k\in\mathbb{Z}_0^2} \frac{\log|k|}{|k|^{2+2r}} < \infty.
\end{aligned}
$$

(3.6)

In the final calculation we have used Proposition A.1 in the Appendix. We need these detailed calculations in order to obtain the estimate in (3.6). Indeed, the present literature (see the references given before Proposition 3.3) deals with the components $B_k$'s, without taking too much care on how the value of $\int |B_k|^2 \, d\mu_\nu$ depends on the index $k$. But this is important for the estimate of the "vector" $B = \sum_k B_k e_k$.

Let us come back to the question of the definition on $B(u)$ in $\mathcal{H}_2^{-r-1}$. By similar calculation as above, one shows that

$$
\|B^N(u) - B(u)\|_{\mathcal{H}_2^{-r-1}} \to 0 \qquad \text{in } \mathcal{L}^2(\mu_\nu), \text{ as } N\to\infty;
$$

therefore, for some subsequence we have that

$$
\|B^N(u) - B(u)\|_{\mathcal{H}_2^{-r-1}} \to 0 \qquad \text{for } \mu_\nu\text{-a.e. } u,
$$

which shows that $B(u)$ is indeed in $\mathcal{H}_2^{-r-1}$. Since $\mu_\nu$ is Gaussian [and bearing in mind (3.3)], similar calculations hold for any even exponent $\rho$ and then by Hölder inequality for any $1 \leq \rho < \infty$.  $\square$

REMARK 3.2.  According to the latter result, the nonlinear term $B(u)$ is defined for $\mu_\nu$-a.e. $u$. Since $\mu_\nu(\mathcal{H}_2^0) = 0$ but $\mu_\nu(\mathcal{H}_q^{-r}) = 1$ ($r > 0$, $1 < q < \infty$), the elements $u$ for which the nonlinear term $B(u)$ exists are (nonregular) distributions. Da Prato and Debussche [9] explain that $B(u) \in \mathcal{L}^\rho(\mu_\nu; \mathcal{H}_2^{-r-1})$ for $1 \leq \rho < \infty$, $r > 0$, as follows. Denote by $:u\otimes u:$ the renormalized square (Wick square), defined as $:u \otimes u: := u \otimes u - \mathbb{E}_{\mu_\nu}(u \otimes u)$ (see, e.g., [12]). Consider the finite-dimensional approximations $u_N = \sum_{|k| \le N} u_k e_k$; one has that $\sup_N \mathbb{E}_{\mu_\nu} \| :u_N \otimes u_N: \|^\rho_{\mathcal{H}_2^{-r}} < \infty$. Notice that $\nabla \cdot (:u_N \otimes u_N:) = \nabla \cdot (u_N \otimes u_N - \mathbb{E}_{\mu_\nu}(u_N \otimes u_N)) = \nabla \cdot (u_N \otimes u_N)$. Hence $B(u_N) = \Pi[\nabla \cdot (:u_N \otimes u_N:)]$ and in the limit $B(u) = \Pi[\nabla \cdot (:u \otimes u:)]$ is well defined, that is $B(u) \in \mathcal{L}^\rho(\mu_\nu; \mathcal{H}_2^{-r-1})$.

Finally, let us recall the main properties of the components $B_k$. (For the proof, see [2, 4, 8]. As noticed above, the proof consists in getting uniform estimates for the sequence of finite approximations $B_k^N$.)

PROPOSITION 3.3. *For any $k \in \mathbb{Z}_0^2$,*

(3.7) $$\partial_k B_k = 0,$$

(3.8) $$\overline{B_k} = B_{-k},$$

(3.9) $$B_k \in \mathcal{L}^p(\mu_\nu) \qquad \text{for any } 1 \le p < \infty.$$

*Each component $B_k$ is the $\mathcal{L}^p(\mu_\nu)$-limit (as $N \to \infty$) of the Galerkin approximations*

$$B_k^N(u) = \sum_{\substack{h \\ 0 < |h|, |k-h|, |k| \le N}} c_{h,k} u_h u_{k-h}, \qquad k \in \mathbb{Z}_0^2, \ N \in \mathbb{N},$$

*for which one has the conservation of the enstrophy, that is,*

$$\sum_{\substack{k \\ 0 < |k| \le N}} B_k^N(u)|k|^2 \overline{u_k} = 0, \qquad N \in \mathbb{N}.$$

## 4. Pathwise uniqueness.
First, we recall the result in [9]. These authors show that there exists a unique (strong) solution $u^x$ of (2.4) for $\mu_\nu$-a.e. $x \in \mathcal{B}_{pq}^{-s}$ (if the parameters satisfy: $0 < s < \frac{2}{p}, 2 < p = q < \infty, s + \frac{2}{p} < 1$; therefore the set of initial data has $\mu_\nu$-measure equal to 1), such that

(4.1) $$u^x - z \in C([0, \infty); \mathcal{B}_{pq}^{-s}) \cap L_{\text{loc}}^\beta([0, \infty); \mathcal{B}_{pq}^\alpha)$$

$\mathbb{P}$-a.s., where $\alpha$ and $\beta$ are suitable parameters and $z$ is the stationary solution of the Ornstein–Uhlenbeck equation (see Section 5 in [9]). Since $z \in C([0, \infty); \mathcal{B}_{pq}^{-s})$ ($\mathbb{P}$-a.s.), then the regularity of $u^x$ (which is the important unknown variable) is obtained by merging together the regularity of $u^x - z$ and of $z$. Therefore the result of [9] states that there exists a process $u^x$ such that

$$u^x \in C([0, \infty); \mathcal{B}_{pq}^{-s}) \qquad \mathbb{P}\text{-a.s.;}$$

moreover, only one of the processes in the space $C([0, \infty); \mathcal{B}_{pq}^{-s})$ satisfies the further condition (4.1).

Finally, this solution $u = \{u^x\}_x$, as well as any other solution obtained as the limit of a subsequence of Galerkin approximations (taking the limit as done in [9]), has invariant measure $\mu_\nu$, in the sense that

$$(4.2) \qquad \int \mathbb{E} f(u^x(t)) \, d\mu_\nu(x) = \int f(x) \, d\mu_\nu(x) \qquad \forall f \in \mathcal{L}^1(\mu_\nu), \ t \geq 0.$$

The fact that $\mu_\nu$ is invariant for the Galerkin approximations $u_N$ is an important tool in the proof of the existence (in the spaces considered in [9] as well as in those considered in [1]). Moreover, any solution $u$, obtained as the limit of a subsequence of Galerkin approximations, has $\mu_\nu$ as invariant measure. Since in both articles [1, 9] the limit of a subsequence is considered and not that of the whole sequence $\{u_N\}$, it is natural to ask about uniqueness of this limit obtained from any subsequence of Galerkin approximations.

We intend, however, to show pathwise uniqueness for solutions with paths in $C([0, \infty); \mathcal{B}_{pq}^{-s})$ and with invariant measure $\mu_\nu$, without the additional requirement (4.1) on $u^x - z$. The invariance of the measure $\mu_\nu$ is used in order to deal with the nonlinear term $B(u)$.

From now on, we consider as state space any Besov space $\mathcal{B}_{pq}^{-s}$ of full measure $\mu_\nu$.

For $\mu_\nu$-a.e. $x \in \mathcal{B}_{pq}^{-s}$ [i.e., $x \in S' \cap \mathcal{B}_{pq}^{-s}$, with $\mu_\nu(S') = 1$], let $u^x$ be a process solving (2.4) such that (4.2) holds and $\mathbb{P}$-a.e. path

$$(4.3) \qquad\qquad\qquad u^x \in C([0, \infty); \mathcal{B}_{pq}^{-s}).$$

In particular, from the invariance formula (4.2) with $f(x) = \|B(x)\|_{\mathcal{H}_2^{-r-1}}^\rho$, one obtains that

$$(4.4) \qquad \int \int_0^T \mathbb{E}\|B(u^x(t))\|_{\mathcal{H}_2^{-r-1}}^\rho \, dt \, d\mu_\nu(x) = T \int \|B(x)\|_{\mathcal{H}_2^{-r-1}}^\rho \, d\mu_\nu(x).$$

[Actually, this holds if $u^x \in C([0, \infty); S)$, $x \in S$, for any $S \subset \mathcal{U}'$ with $\mu_\nu(S) = 1$.]

Because of (3.5), the quantity on the right-hand side is finite for any finite time $T$ and any $r > 0$, $1 \leq \rho < \infty$. Fix now these parameters. From the left-hand side of (4.4), we obtain that there exists a subset $S'' \subset (S' \cap \mathcal{B}_{pq}^{-s}) \subset \mathcal{U}'$ with $\mu_\nu(S'') = 1$ such that

$$\forall x \in S'' \qquad \mathbb{E} \int_0^T \|B(u^x(t))\|_{\mathcal{H}_2^{-r-1}}^\rho \, dt < \infty,$$

and therefore

$$\forall x \in S'' \ \exists \Omega^x \subset \Omega, \ \mathbb{P}(\Omega^x) = 1 \colon \int_0^T \|B(u^x(t, \omega))\|_{\mathcal{H}_2^{-r-1}}^\rho \, dt < \infty \qquad \forall \omega \in \Omega^x.$$

We repeat this procedure for a countable choice of the parameters ($\rho$, $T = 1, 2, \ldots$; $r = \frac{1}{2}, \frac{1}{3}, \ldots$) and use interpolation results for all positive real numbers $r \neq \frac{1}{2}, \frac{1}{3}, \ldots$. Then we obtain that, for $\mu_\nu(S'') = 1$,

$$\forall x \in S'' \,\, \exists \Omega^x \subset \Omega, \,\, \mathbb{P}(\Omega^x) = 1 : \int_0^T \|B(u^x(t,\omega))\|_{\mathcal{H}_2^{-r-1}}^\rho \, dt < \infty$$

(4.5)

$$\forall \omega \in \Omega^x, \,\, T > 0, \,\, r > 0, \,\, 1 \le \rho < \infty.$$

Hence, given $x \in S''$, the solution $u^x$ enjoys ($\mathbb{P}$-a.s.) the property

(4.6) $$\int_0^T \|B(u^x(t))\|_{\mathcal{H}_2^{-r-1}}^\rho \, dt < \infty \qquad \forall T > 0, \,\, r > 0, \,\, 1 \le \rho < \infty.$$

Let $\tilde{u}^x$ be any other process defined on the same probability space ($\Omega$, $\mathcal{F}$, $\{\mathcal{F}_t\}$, $\mathbb{P}$), with the same properties given above for $u^x$ and solving (2.4) with the same $\{\mathcal{F}_t\}$-Wiener process as for $u^x$. Define the difference $v^x = u^x - \tilde{u}^x$; then $v^x \in C([0,\infty); \mathcal{B}_{pq}^{-s})$. From now on we drop the dependence on $x$. $v$ satisfies the equation

(4.7) $$\frac{d}{dt}v(t) + Av(t) = -B(u(t)) + B(\tilde{u}(t)), \qquad t > 0,$$
$$v(0) = 0.$$

Bearing in mind the regularizing effect of the Stokes operator $A$, something more can be proven. More precisely, (4.6) grants that the right-hand side of the first equation in (4.7) belongs to the space $L_{\text{loc}}^\rho(0,\infty; \mathcal{H}_2^{-r-1})$ for any $1 \le \rho < \infty$, $r > 0$. By Proposition A.2 in the Appendix, one has that

(4.8) $$v \in L_{\text{loc}}^\rho(0,\infty; \mathcal{H}_2^{-r+1}) \cap C([0,\infty); \mathcal{B}_{2\rho}^{-r+1-2/\rho}).$$

This holds for any $r > 0$, $1 < \rho < \infty$. Hence we have proven that any solution $v$ to (4.7) must belong to the functional space $\Sigma := \bigcap_{1 < \rho < \infty, r > 0} \Sigma_{\rho,r}$, where $\Sigma_{\rho,r} := L_{\text{loc}}^\rho(0,\infty; \mathcal{H}_2^{-r+1}) \cap C([0,\infty); \mathcal{B}_{2\rho}^{-r+1-2/\rho})$. Let us point out that, for $2 \le p \le \rho \le q$, we have $\mathcal{B}_{2\rho}^{-r+1-2/\rho} \subseteq \mathcal{B}_{\rho\rho}^{-r-2/\rho+2/p} \subseteq \mathcal{B}_{p\rho}^{-r} \subseteq \mathcal{B}_{pp}^{-r}$, and for $r \le s$, we have $\mathcal{B}_{pq}^{-r} \subseteq \mathcal{B}_{pq}^{-s}$; therefore, $\mathcal{B}_{2\rho}^{-r+1-2/\rho} \subseteq \mathcal{B}_{pq}^{-s}$. Thus the regularity specified in (4.8) is stronger than the regularity $v \in C([0,\infty); \mathcal{B}_{pq}^{-s})$ given by the definition of $v$ itself, as $v = u - \tilde{u}$.

REMARK 4.1. The regularizing effect of the Stokes operator is not enough to obtain more regularity in the stochastic equation (2.4), because of the presence of the cylindric noise $dw$. This is already evident for the stochastic Stokes equation, that is, the equation obtained from (2.4) by neglecting the nonlinear operator $B$ (see, e.g., [9] for the optimal regularity of the stochastic Stokes equation, where it is shown that the solution $z$ of the stochastic Stokes equation does take values in distribution spaces). Therefore $u$, as well as $z$, are expected to have paths in $C([0,\infty); \mathcal{B}_{pq}^{-s})$.

Bearing in mind the bilinearity of the operator $B$, the equation for $v$ can be written in the following form:

$$(4.9) \qquad \frac{d}{dt}v(t) + Av(t) + B(u(t), v(t)) + B(v(t), \tilde{u}(t)) = 0, \qquad t > 0,$$
$$v(0) = 0.$$

REMARK 4.2. Actually, so far the equivalence between (4.7) and (4.9) holds only heuristically. Of course, for the rigorous equivalence of this equality it is necessary that $B(u, v) + B(v, \tilde{u})$ is meaningful. We shall see in the proof of the next theorem that this is indeed the case, because $v$ is more regular than $u$ and $\tilde{u}$, as already shown in (4.8).

More precisely, for $dt$-a.e. $t \in [0, T]$, for the $N$-finite-dimensional approximations we have

$$B(u_N^x(t), u_N^x(t)) - B(\tilde{u}_N^x(t), \tilde{u}_N^x(t)) = B(u_N^x(t), v_N^x(t)) + B(v_N^x(t), \tilde{u}_N^x(t)).$$

The left-hand side converges to $B(u^x(t)) - B(\tilde{u}^x(t))$; indeed, proceeding as in Section 3, we prove that $B(u_N^x) \to B(u^x)$ in $\mathcal{L}^1(\mu_\nu; L_{\mathrm{loc}}^\rho(0, \infty; \mathcal{H}_2^{-r-1}))$, as $N \to \infty$ (for $r > 0$ and $1 \le \rho < \infty$), and hence, for $\mu_\nu \times dt$-a.e. $(x, t)$, some subsequence of $B(u_N^x(t))$ converges to $B(u^x(t))$ in $\mathcal{H}_2^{-r-1}$, as $N \to \infty$. The same holds for $\tilde{u}$.

The right-hand side has a limit, thanks to the regularity of $v$. In particular, under the assumptions (4.10), for fixed $x \in \mathcal{B}_{pq}^{-s}$ and for $dt$-a.e. $t$, $u^x(t) \in \mathcal{B}_{pq}^{-s}$ and $v^x(t) \in \mathcal{B}_{pq}^a$. Then Chemin's estimate [6] on the bilinear operator $B$, as in the proof of the next theorem, gives that the expression $B(u^x(t), v^x(t))$ exists and $B(u_N^x(t), v_N^x(t)) \to B(u^x(t), v^x(t))$ in $\mathcal{B}_{pq}^{-s+a-2/p-1}$, as $N \to \infty$.

The function $v \equiv 0$ is a solution to (4.9). We are going to prove that this is the only solution of (4.9) in the class $\Sigma$.

To prove this, we first show that, given $u, \tilde{u} \in C([0, \infty); \mathcal{B}_{pq}^{-s})$, under the assumptions (4.10), there exists a unique solution $v$ to the problem (4.9) into a class less regular than $\Sigma$. This is proven in Theorem 4.1. From this, uniqueness in the smaller class $\Sigma$ immediately follows. This concludes our proof that the unique solution for (4.7) is $v \equiv 0$. What remains to be proven is therefore the following.

THEOREM 4.1. *Let real numbers $s, a$ be given as well as $1 < \alpha, p, q < \infty$ satisfying the following conditions:*

$$0 < s < a,$$
$$a < \frac{2}{p},$$

(4.10)

$$\frac{1}{2}\left(s+\frac{2}{p}+1\right) < 1,$$

$$\frac{1}{2}\left(-a+\frac{2}{p}+1\right)\frac{\alpha}{\alpha-1} < 1.$$

*Then, for any $u, \tilde{u} \in C([0,\infty); \mathcal{B}_{pq}^{-s})$, there exists a unique $v \in \mathcal{V} := C([0,\infty); \mathcal{B}_{pq}^{-s}) \cap L_{\mathrm{loc}}^{\alpha}(0,\infty; \mathcal{B}_{pq}^{a})$ solution to the following problem:*

(4.11)
$$\frac{d}{dt}v(t) + Av(t) + B(u(t), v(t)) + B(v(t), \tilde{u}(t)) = 0, \qquad t > 0,$$
$$v(0) = 0.$$

*In particular, if $v$ satisfies (4.11), then $v(t) = 0$ for $t \geq 0$.*

PROOF. To begin with, we fix any finite time interval $[0, T]$. We consider the solution to (4.11) in the mild form (in the sense of, e.g., [10])

(4.12)
$$v(t) = -\int_0^t e^{-(t-\tau)A}[B(u(\tau), v(\tau)) + B(v(\tau), \tilde{u}(\tau))]\, d\tau.$$

We want to prove existence and uniqueness of a solution in $\mathcal{V}_T := C([0, T]; \mathcal{B}_{pq}^{-s}) \cap L^{\alpha}(0, T; \mathcal{B}_{pq}^{a})$ by a fixed point theorem, as in [9]. We consider the norm $\|v\|_{\mathcal{V}_T} = \|v\|_{C([0,T]; \mathcal{B}_{pq}^{-s})} + \|v\|_{L^{\alpha}(0,T; \mathcal{B}_{pq}^{a})}$. We proceed in three steps.

*Step* 1. We begin by estimating the bilinear operator by means of Bony's para-products techniques, as given in [6], Corollary 1.3.1:

(4.13)
$$\|B(u,v)\|_{\mathcal{B}_{pq}^{-s+a-2/p-1}} = \|\nabla \cdot (u \otimes v)\|_{\mathcal{B}_{pq}^{-s+a-2/p-1}}$$
$$\leq \|u \otimes v\|_{\mathcal{B}_{pq}^{-s+a-2/p}}$$
$$\leq c\|u\|_{\mathcal{B}_{pq}^{-s}}\|v\|_{\mathcal{B}_{pq}^{a}},$$

if

(4.14)
$$0 < s < a \quad \text{and} \quad a < \frac{2}{p}.$$

We remark that $B(u, v)$ makes sense, when at least one element belongs to a Besov space of positive order ($v \in \mathcal{B}_{pq}^{a}$ with $a > 0$).

*Step* 2. Let us show that, given $v \in \mathcal{V}_T$, the right-hand side of (4.12) belongs to $\mathcal{V}_T$. By the property of the Stokes operator [basically, the property of the heat operator: $\|e^{-tA}x\|_{\mathcal{B}_{pq}^{a}} \leq ct^{-(a-b)/2}\|x\|_{\mathcal{B}_{pq}^{b}}$ for $t > 0$ and $a \geq b$], the following holds:

$$\left\|\int_0^t e^{-(t-\tau)A}B(u(\tau), v(\tau))\, d\tau\right\|_{\mathcal{B}_{pq}^{a}}$$

$$
\begin{aligned}
(4.15) \qquad & \leq \int_0^t \| e^{-(t-\tau)A} B(u(\tau), v(\tau)) \|_{\mathcal{B}_{pq}^a} \, d\tau \\
& \leq c \int_0^t \frac{1}{(t-\tau)^{(s+2/p+1)/2}} \| B(u(\tau), v(\tau)) \|_{\mathcal{B}_{pq}^{-s+a-2/p-1}} \, d\tau \\
& \leq c \| u \|_{C([0,T]; \mathcal{B}_{pq}^{-s})} \int_0^t \frac{1}{(t-\tau)^{(s+2/p+1)/2}} \| v(\tau) \|_{\mathcal{B}_{pq}^a} \, d\tau
\end{aligned}
$$

(denoting different constants by the same symbol $c$).

We now estimate the convolution integral by Young's inequality. Thus

$$
\begin{aligned}
(4.16) \qquad & \left\| \int_0^t e^{-(t-\tau)A} B(u(\tau), v(\tau)) \, d\tau \right\|_{L^\alpha(0,T; \mathcal{B}_{pq}^a)} \\
& \leq C_1 T^{(-s-2/p+1)/2} \| u \|_{C([0,T]; \mathcal{B}_{pq}^{-s})} \| v \|_{L^\alpha(0,T; \mathcal{B}_{pq}^a)},
\end{aligned}
$$

if

$$
(4.17) \qquad \frac{1}{2} \left( s + \frac{2}{p} + 1 \right) < 1.
$$

In the same way, we check the estimate in $C([0,T]; \mathcal{B}_{pq}^{-s})$. First

$$
\begin{aligned}
(4.18) \qquad & \left\| \int_0^t e^{-(t-\tau)A} B(u(\tau), v(\tau)) \, d\tau \right\|_{\mathcal{B}_{pq}^{-s}} \\
& \leq c \int_0^t \frac{1}{(t-\tau)^{(-a+2/p+1)/2}} \| B(u(\tau), v(\tau)) \|_{\mathcal{B}_{pq}^{-s+a-2/p-1}} \, d\tau \\
& \leq c \| u \|_{C([0,T]; \mathcal{B}_{pq}^{-s})} \int_0^t \frac{1}{(t-\tau)^{(-a+2/p+1)/2}} \| v(\tau) \|_{\mathcal{B}_{pq}^a} \, d\tau.
\end{aligned}
$$

Again Young's inequality allows us to conclude that the latter expression is finite if

$$
(4.19) \qquad \frac{1}{2} \left( -a + \frac{2}{p} + 1 \right) \frac{\alpha}{\alpha - 1} < 1,
$$

and moreover,

$$
\begin{aligned}
(4.20) \qquad & \left\| \int_0^t e^{-(t-\tau)A} B(u(\tau), v(\tau)) \, d\tau \right\|_{L^\infty([0,T]; \mathcal{B}_{pq}^{-s})} \\
& \leq C_2 T^{((a-2/p-1)/2)(\alpha/(\alpha-1))+1} \| u \|_{C([0,T]; \mathcal{B}_{pq}^{-s})} \| v \|_{L^\alpha(0,T; \mathcal{B}_{pq}^a)}.
\end{aligned}
$$

We notice that the same computations hold for $B(v, \tilde{u})$.

Hence, if $v \in \mathcal{V}_T$ and (4.14), (4.17), (4.19) hold, then $\int_0^t e^{-(t-\tau)A} B(u(\tau), v(\tau)) \, d\tau \in \mathcal{V}_T$.

*Step* 3. Equation (4.11) is linear in $v$. Hence the estimates (4.16) and (4.20) give that the mapping

$$v \mapsto -\int_0^t e^{-(t-\tau)A}[B(u(\tau), v(\tau)) + B(v(\tau), \tilde{u}(\tau))] \, d\tau$$

is a contraction in $\mathcal{V}_{T^*}$ with $T^* \leq T$ and such that

$$(4.21) \quad T^* < \min\left\{ \left(\frac{1}{2C_1 N_T}\right)^{1/(((a-2/p-1)/2)(\alpha/(\alpha-1))+1)}, \right.$$
$$\left. \left(\frac{1}{2C_2 N_T}\right)^{1/((-s-2/p+1)/2)} \right\},$$

where $N_T = \|u\|_{C([0,T]; \mathcal{B}_{pq}^{-s})} + \|\tilde{u}\|_{C([0,T]; \mathcal{B}_{pq}^{-s})}$. Hence, on the interval $[0, T^*)$, there exists a unique solution $v$ with the regularity specified in $\mathcal{V}$. This is $v(t) = 0$ for $0 \leq t < T^*$. Notice that the amplitude of the time interval for local existence depends only on the $C([0,T]; \mathcal{B}_{pq}^{-s})$-norms of $u$ and $\tilde{u}$; therefore, we can continue in such a way as to cover the time interval $[0, T]$ with a finite number of intervals of amplitude $\frac{3}{4}T^*$.

Since this holds for any finite $T$, the proof is completed. $\quad\square$

REMARK 4.3. Since $s > 0$, the third condition on (4.10) imposes that $p > 2$. This is the reason for working in Besov spaces, instead of the usual Hilbert spaces.

Choose now the parameters of Theorem 4.1 to be $p = q = \alpha = 3$, $s = \frac{1}{6}$, $a = \frac{1}{2}$. In this way, bearing in mind Proposition 3.1, we have fixed a set $\mathcal{B}_{pq}^{-s}$ of initial data such that $\mu_\nu(\mathcal{B}_{pq}^{-s}) = 1$ (but many other choices are possible); moreover, the assumptions of Theorem 4.1 are satisfied. Choose also the parameters $\rho = 3$, $r = \frac{1}{6}$ for the regularity of (4.8). Finally, by an embedding theorem [see [5], Theorem 6.5.1], we have

$$\mathcal{B}_{2\rho}^{-r+1-2/\rho} \subset \mathcal{B}_{pq}^{-s},$$
$$\mathcal{H}_2^{-r+1} \subset \mathcal{B}_{pq}^a.$$

Hence $\Sigma \subset \mathcal{V}$. And the uniqueness in $\mathcal{V}$ implies the uniqueness in $\Sigma$.

We have therefore proven the following.

THEOREM 4.2. *Pathwise uniqueness of the solutions to the stochastic Navier–Stokes equation with space-time Gaussian white noise* (2.4), *for which $\mu_\nu$ is an invariant measure, holds in the following precise sense: there exists a set $S \subset \mathcal{U}'$ with $\mu_\nu(S) = 1$ such that for, $\mu_\nu$-a.e. $x \in S$, the $C([0,\infty); S)$-valued paths of any two solutions of* (2.4), *defined on the same probability space with the same Wiener process and having invariant measure $\mu_\nu$, coincide $\mathbb{P}$-a.s.*

## APPENDIX

In this appendix, two results used in the previous proofs are presented. We begin with the estimate on the sum of the series $\sum_{h \in \mathbb{Z}_0^2, h \neq k} \frac{1}{|h|^2 |k-h|^2}$ which is (absolutely) convergent for each $k \in \mathbb{Z}_0^2$. It is enough to perform the calculation for the integral

$$\int_{\mathbb{R}^2 \setminus (C_0 \cup C_k)} \frac{1}{(x^2 + y^2)([x - k_1]^2 + [y - k_2]^2)} \, dx \, dy$$

with $C_h = \{(x, y) \in \mathbb{R}^2 : [x - h_1]^2 + [y - h_2]^2 \leq 1\}$, given $h = (h_1, h_2) \in \mathbb{Z}^2$. By a rotation around the origin bringing the point $k$ into the semipositive $x$-axis (so $C_k$ is $C_{(|k|,0)}$), the integral can be written as

$$(A.1) \qquad \int_{\mathbb{R}^2 \setminus (C_0 \cup C_k)} \frac{1}{(x^2 + y^2)([x - |k|]^2 + y^2)} \, dx \, dy.$$

We state the following.

PROPOSITION A.1.  *There exists a positive constant $c$ such that*

$$\int_{\mathbb{R}^2 \setminus (C_0 \cup C_k)} \frac{1}{(x^2 + y^2)([x - |k|]^2 + y^2)} \, dx \, dy \leq c \frac{\log |k|}{|k|^2} \qquad \forall k \in \mathbb{Z}_0^2, \ |k| \geq 2.$$

PROOF.  The proof is based on elementary calculations. We show the main steps. First, we note that the integrand function can be written as the sum of four terms:

$$
\begin{aligned}
(A.2) \qquad & \frac{1}{(x^2 + y^2)([x - |k|]^2 + y^2)} \\
& = \frac{2x}{|k|(x^2 + y^2)(|k|^2 + 4y^2)} - \frac{2[x - |k|]}{|k|([x - |k|]^2 + y^2)(|k|^2 + 4y^2)} \\
& \quad + \frac{1}{(x^2 + y^2)(|k|^2 + 4y^2)} + \frac{1}{([x - |k|]^2 + y^2)(|k|^2 + 4y^2)}.
\end{aligned}
$$

For the integral of the second addendum, one has

$$\int_{\mathbb{R}^2 \setminus (C_0 \cup C_k)} \frac{2[x - |k|]}{|k|([x - |k|]^2 + y^2)(|k|^2 + 4y^2)} \, dx \, dy$$

$$= \int_{\mathbb{R}^2 \setminus (C_{-k} \cup C_0)} \frac{2x}{|k|(x^2 + y^2)(|k|^2 + 4y^2)} \, dx \, dy$$

by a change of variable. Therefore the integrals of the first two addenda on the right-hand side of (A.2) partly cancel each other, and what is left are two integrals on small balls:

$$\int_{C_{-k}} \frac{2x}{|k|(x^2 + y^2)(|k|^2 + 4y^2)} \, dx \, dy - \int_{C_k} \frac{2x}{|k|(x^2 + y^2)(|k|^2 + 4y^2)} \, dx \, dy.$$

This quantity vanishes, by symmetry.

Hence, the only contribution to the integral (A.1) comes from the last two addenda in (A.2). We have

$$\int_{\mathbb{R}^2 \setminus (C_0 \cup C_k)} \frac{dx\, dy}{(x^2 + y^2)(|k|^2 + 4y^2)} \leq \int_{\mathbb{R}^2 \setminus C_0} \frac{dx\, dy}{(x^2 + y^2)(|k|^2 + 4y^2)}$$

and

$$\int_{\mathbb{R}^2 \setminus (C_0 \cup C_k)} \frac{dx\, dy}{([x - |k|]^2 + y^2)(|k|^2 + 4y^2)}$$

$$= \int_{\mathbb{R}^2 \setminus (C_{-k} \cup C_0)} \frac{dx\, dy}{(x^2 + y^2)(|k|^2 + 4y^2)}$$

$$\leq \int_{\mathbb{R}^2 \setminus C_0} \frac{dx\, dy}{(x^2 + y^2)(|k|^2 + 4y^2)}.$$

It remains to calculate this latter integral. We proceed as follows. Let $Q_0$ be the rectangle $\{(x, y) \in \mathbb{R}^2 : |x| \leq \frac{1}{\sqrt{2}}, 0 < y \leq \frac{1}{\sqrt{2}}\}$. Then

$$\int_{\mathbb{R}^2 \setminus C_0} \frac{dx\, dy}{(x^2 + y^2)(|k|^2 + 4y^2)}$$

$$\leq 2 \int_{\mathbb{R}^2 \setminus Q_0} \frac{dx\, dy}{(x^2 + y^2)(|k|^2 + 4y^2)}$$

$$= 2 \int_0^{1/\sqrt{2}} \frac{dy}{|k|^2 + 4y^2} \int_{|x| \geq 1/\sqrt{2}} \frac{dx}{x^2 + y^2} + 2 \int_{1/\sqrt{2}}^{\infty} \frac{dy}{|k|^2 + 4y^2} \int_{\mathbb{R}} \frac{dx}{x^2 + y^2}.$$

Let us estimate these two integrals. For the first, we have

$$\int_0^{1/\sqrt{2}} \frac{dy}{|k|^2 + 4y^2} \int_{|x| \geq 1/\sqrt{2}} \frac{dx}{x^2 + y^2}$$

$$\leq \int_0^{1/\sqrt{2}} \frac{dy}{|k|^2 + 4y^2} \int_{|x| \geq 1/\sqrt{2}} \frac{dx}{x^2}$$

$$\leq \int_0^{1/\sqrt{2}} \frac{dy}{|k|^2} 2\sqrt{2}$$

$$= \frac{2}{|k|^2}.$$

For the second,

$$\int_{1/\sqrt{2}}^{\infty} \frac{dy}{|k|^2 + 4y^2} \int_{\mathbb{R}} \frac{dx}{x^2 + y^2} = 2 \int_{1/\sqrt{2}}^{\infty} \frac{1}{|k|^2 + 4y^2} \frac{\pi}{y}\, dy$$

$$= 2\pi \int_{1/\sqrt{2}}^{|k|} \frac{dy}{(|k|^2 + 4y^2)y} + 2\pi \int_{|k|}^{\infty} \frac{dy}{(|k|^2 + 4y^2)y}$$

$$\leq \frac{2\pi}{|k|^2+2} \int_{1/\sqrt{2}}^{|k|} \frac{dy}{y} + \frac{2\pi}{|k|} \int_{|k|}^{\infty} \frac{dy}{|k|^2+4y^2}$$

$$\leq 2\pi \frac{\log(\sqrt{2}|k|)}{|k|^2+2} + \frac{2\pi}{|k|} \frac{1}{2|k|} \frac{\pi}{2}.$$

Summing up all the estimates, the proof is completed.  □

The second result concerns regularity for parabolic equations.

PROPOSITION A.2.  *Let $T \in (0,\infty]$, $1 < \rho < \infty$, and $\sigma \in \mathbb{R}$. Let $A$ be the Stokes operator described in Section 2.*
*For any $f \in L^\rho(0,T;\mathcal{H}_2^\sigma)$, the Cauchy problem*

$$\frac{d}{dt}X(t) + AX(t) = f(t), \qquad t \in (0,T],$$

$$X(0) = 0,$$

*has a unique solution $X \in \mathcal{W}^{1,\rho}(0,T) \equiv \{X \in L^\rho(0,T;\mathcal{H}_2^{\sigma+2}) : \frac{d}{dt}X \in L^\rho(0,T; \mathcal{H}_2^\sigma)\}$. Moreover, the solution depends continuously on the data in the sense that there exists a constant $c_{\rho,\sigma}$ such that*

$$\left( \int_0^T \left[ \|X(t)\|_{\mathcal{H}_2^{\sigma+2}}^\rho + \left\| \frac{d}{dt}X(t) \right\|_{\mathcal{H}_2^\sigma}^\rho \right] dt \right)^{1/\rho} \leq \left( c_{\rho,\sigma} \int_0^T \|f(t)\|_{\mathcal{H}_2^\sigma}^\rho \, dt \right)^{1/\rho}.$$

*Finally, $X \in C_b([0,T]; \mathcal{B}_{2\rho}^{\sigma+2-2/\rho})$.*

PROOF.  The Stokes operator $A$ is a positive self-adjoint operator in $\mathcal{H}_2^\sigma$ with domain $\mathcal{H}_2^{\sigma+2}$, and it generates an analytic semigroup in $\mathcal{H}_2^\sigma$. Then the first part of the proposition is obtained applying Theorem 3.2 in [11]. Moreover, by interpolation we get that the space $\mathcal{W}^{1,\rho}(0,T)$ is continuously embedded in the space $C_b([0,T]; \mathcal{B}_{2\rho}^{\sigma+2-2/\rho})$; that is, there exists a positive constant $c$ such that

$$\|X\|_{C_b([0,T];\mathcal{B}_{2\rho}^{\sigma+2-2/\rho})} \leq c \|X\|_{\mathcal{W}^{1,\rho}(0,T)}. \qquad \qquad \square$$

**Acknowledgments.**  We would like to thank Zdzisław Brzeźniak for very stimulating discussions. The second author thanks the Institut für Angewandte Mathematik of Bonn University for the warm hospitality.

INSTITUT FÜR ANGEWANDTE
MATHEMATIK
UNIVERSITÄT BONN
WEGELERSTRASSE 6
D-53115 BONN
GERMANY
E-MAIL: albeverio@uni-bonn.de

DIPARTIMENTO DI MATEMATICA
UNIVERSITÀ DI PAVIA
VIA FERRATA 1
I-27100 PAVIA
ITALY
E-MAIL: ferrario@dimat.unipv.it